\documentclass[a4paper,12pt]{article}

\usepackage{amsmath,amsfonts,amssymb,bm,amsthm}
\usepackage{mathrsfs}
\usepackage{cases}
\usepackage{cite}
\usepackage{indentfirst}
\usepackage{fancyhdr}
\usepackage{geometry}
\usepackage[usenames]{color}
\usepackage[colorlinks=true,
linkcolor=webblack,
filecolor=webblack,
citecolor=webblack]{hyperref}

\definecolor{webblack}{rgb}{0,.5,0}
\definecolor{webblack}{rgb}{.6,0,0}

\usepackage{color}

\newtheorem{theorem}{Theorem}

\newtheorem{definition}{Definition}

\newtheorem{corollary}{Corollary}

\allowdisplaybreaks

\begin{document}

\title{\textbf{Geometrodynamics based on geodesic equation with Cartan structural equation on Riemannian manifolds }}
\vspace{.3 cm}
\author{Gen Wang\thanks{\tt
Wg1991dream@163.com
}}
\date{{\em \small School of Mathematical Sciences, Xiamen University,\\
     Xiamen, 361005, P.R.China.}}

\maketitle
\abstract{Motivated by the geospin matrix $W$ as a new variable in Riemannian geometry.  Then, we use four real dynamical variables $\left\{ a,\alpha ,v,W \right\}$  to show the dynamical essence of Cartan structural equation, we obtain the geometrodynamics on Riemannian manifolds that can be expressed below
\begin{align}
  & \Theta /d{{t}^{2}}=a-v\wedge W \notag\\
 & d\Theta /d{{t}^{3}}=v\wedge \alpha -a\wedge W\notag\\
 & \Omega /d{{t}^{2}}=\alpha -W\wedge W \notag\\
 & d\Omega /d{{t}^{3}}=W\wedge \alpha -\alpha \wedge W\notag
\end{align}that is valid more generally for any connection in a principal bundle. We can see that the first equation explains the geodesic equation very well, the second formula means the first Bianchi identity, while the last equation reveals the dynamical nature of the second Bianchi identity. It implies that Cartan structural equations on Riemannian manifolds can be rewritten in a real geometrodynamical form based on the geospin matrix.

}


\newpage
\tableofcontents

\section{Introduction}
\label{sec:intro}
In Riemannian manifolds \cite{1}, the most important concept is the metric,
with the help of the metric, the entity of geodesic line can be introduced naturally. In a narrow sense, it is the shortest line between two points and the extension of the concept of straight line segment in a flat space, the  geodesic line fits the general space-- the curved spacetime,  the mathematical expression describing geodesic is called geodesic equation which can be derived by many methods.   Geodesic equation is a system of second order differential equations. Because the unknown functions and their first derivatives are coupled in each equation, the solution is not simple \cite{2,3}. The most prominent application is in general relativity \cite{4,5,6,7} that a geodesic generalizes the notion of a straight line to curved spacetime. Importantly, the world line of a particle free from all external, non-gravitational forces is a particular type of geodesic. In other words, a freely moving or falling particle always moves along a geodesic.
Gravity \cite{8,9,10,11,12} can be regarded as not a force but a consequence of a curved spacetime geometry where the source of curvature is the stress-energy tensor. Thus, for example, the path of a planet orbiting a star is the projection of a geodesic of the curved four-dimensional spacetime geometry around the star onto three-dimensional space. Once the spacetime metric is determined by Einstein field equation, the geodesic equation can be solved to obtain the trajectory of a free particle. Einstein believed that gravity is not a force, but a result of the bending of four-dimensional space by mass and energy.  Above all,  when light passes through the vicinity of a massive object, it will deflect due to its gravitation or enter the curved space near the object, which is called geodesic effect. The geodetic effect represents the effect of the curvature of spacetime, predicted by general relativity, on a vector carried along with an orbiting body.  The term geodetic effect has two slightly different meanings as the moving body may be spinning or non-spinning. Non-spinning bodies move in geodesics, whereas spinning bodies move in slightly different orbits.   One can attempt to explain the geodetic effect by using the  geospin variables \cite{13}.

\subsection{Overviews of geodesic equation}
In this part, we firstly give a simple overviews of geodesic equation for the convenience of later discussion.
In Riemannian geometry, the geodesic \cite{1,3,5,6} as a generalization of the notion of a straight line to a more general setting is a curve representing in some sense the shortest path between two points in a Riemannian manifold.  In a Riemannian manifold or submanifold geodesics are characterised by the property of having vanishing geodesic curvature. More generally, in the presence of an affine connection, a geodesic is defined to be a curve whose tangent vectors remain parallel if they are transported along it. Applying this to the Levi-Civita connection of a Riemannian metric recovers the previous notion. The geodesic equation is of the form given by \cite{14,15,16}
\[\frac{{{d}^{2}}{{x}^{k}}}{d{{t}^{2}}}+\Gamma _{ij}^{k}\frac{d{{x}^{j}}}{dt}\frac{d{{x}^{i}}}{dt}=0\]

The quantity on the left hand side of this equation is the acceleration of a particle, so this equation is analogous to Newton's laws of motion, which likewise provide formulae for the acceleration of a particle. This equation of motion employs the Einstein notation, meaning that repeated indices are summed. The Christoffel symbols are functions of the four spacetime coordinates and so are independent of the velocity or acceleration or other characteristics of a test particle whose motion is described by the geodesic equation. Physicist Steven Weinberg has presented a derivation of the geodesic equation of motion directly from the equivalence principle.

By noticing that the geodesic equation is a second-order ODE which is an ordinary differential equation for the coordinates. It has a unique solution, given an initial position and an initial velocity. Therefore, from the point of view of classical mechanics, geodesics can be thought of as trajectories of free particles in a manifold. Indeed, the equation means that the acceleration vector of the curve has no components in the direction of the surface. So, the motion is completely determined by the bending of the surface. This is also the idea of general relativity where particles move on geodesics and the bending is caused by the gravity. Geodesic should be the shortest path of pure geometry without considering any forces. It also can be written in the form \cite{13}
\begin{equation}
  \left\{ \begin{matrix}
   \frac{d{{x}^{k}}}{dt}={{v}^{k}}\begin{matrix}
   {} &  \\
\end{matrix}  \notag\\
   \frac{d{{v}^{k}}}{dt}+W_{j}^{k}{{v}^{j}}=0 \notag\\
\end{matrix} \right.\notag
\end{equation}
where $W_{j}^{k}$ are geospin variables (see \cite{13} for details).
This is an ordinary differential equation for the coordinates. It has a unique solution, given an initial position and an initial velocity.
By the way, in metric geometry, a geodesic is a curve which is everywhere locally a distance minimizer.
Therefore, from the point of view of classical mechanics, geodesics can be thought of as trajectories $\gamma$ of free particles in a manifold. Indeed, the equation \cite{1,13,15,16}$$\nabla _{\dot {\gamma }}{\dot {\gamma }}=0$$ means that the acceleration vector of the curve has no components in the direction of the surface.  So, the motion is completely determined by the bending of the surface. This is also the idea of general relativity where particles move on geodesics and the bending is caused by the gravity.    Geodesics are of particular importance in general relativity \cite{5,6,7,8,9,10,11}. Timelike geodesics in general relativity describe the motion of free falling test particles.

\section{Riemannian geometry}
Riemannian geometry \cite{1,15,16} is concerned with the objective entity which has nothing to do with all coordinates.
In Riemannian geometry, the Levi-Civita connection is a specific connection on the tangent bundle of a manifold. More specifically, it is the torsion-free metric connection, i.e., the torsion-free connection on the tangent bundle  preserving a given (pseudo-) Riemannian metric. An affine connection is a geometric object on a smooth manifold which connects nearby tangent spaces, and so permits tangent vector fields to be differentiated as if they were functions on the manifold with values in a fixed vector space.

\begin{definition}\cite{1,14,16}
  Suppose that $\nabla:{{C}^{\infty }}\left( TM \right)\times {{C}^{\infty }}\left( TM \right)\to {{C}^{\infty }}\left( TM \right)$ is an operator on tangent bundle, denote $\nabla\left( X,Y \right)={{\nabla}_{X}}Y$ for all $X,Y\in {{C}^{\infty }}\left( TM \right)$. If for vector fields $X,Y,Z$ and functions $f,g$ on manifold $M$ such that\\
  1, ${{\nabla}_{fX+gY}}Z=f{{\nabla}_{X}}Z+g{{\nabla}_{Y}}Z$.\\
  2, ${{\nabla}_{X}}\left( Y+Z \right)={{\nabla}_{X}}Y+{{\nabla}_{X}}Z$.\\
  3, ${{\nabla}_{X}}\left( fY \right)=\left( Xf \right)Y+f{{\nabla}_{X}}Y$.\\
 Then $\nabla$ is said to be a connection on tangent bundle $TM$.

\end{definition}

A smooth Riemannian manifold $(M,g)$ which is a real smooth manifold $M$ equipped with an inner product $g_{p}$ on the tangent space $T_{p}M$ at each point $p$, and it smoothly varies between points that if $X$ and $Y$ are differentiable vector fields on $M$, hence, $p\mapsto g(X,Y)$ is also a smooth function. The family $g_{p}$ of inner products is called a Riemannian metric tensor. Let $M$ be a differentiable manifold of dimension $n$. A Riemannian metric on $M$ is a family of positive definite inner products  \cite{1,16}
$$ g\colon T_{p}M\times T_{p}M\longrightarrow \mathbb{R},\qquad p\in M$$
such that $p\mapsto g(X,Y)$ for all differentiable vector fields $X,Y$ on $M$, then it naturally defines a smooth function $M \rightarrow\mathbb{R}$.  Note that a Riemannian metric $g$ is a symmetric (0,2)-tensor that is positive definite.  Locally, the local coordinates on the manifold $M$ given by ${{x}^{1}},{{x}^{2}},\cdots ,{{x}^{n}}$, the it forms the natural vector fields  $ \left\{{\frac {\partial }{\partial x^{1}}},\dotsc ,{\frac {\partial }{\partial x^{n}}}\right\}$ that 
gives a basis of tangent vectors at each point of $M$. Relative to this coordinate system, the components of the metric tensor are
$$ g_{ij}=g\left( \frac{\partial }{\partial {{x}^{i}}},\frac{\partial }{\partial {{x}^{j}}} \right)$$ at each point $p$.  
Equivalently, the metric tensor can be written in terms of the dual basis $d{{x}^{1}},d{{x}^{2}},\cdots ,d{{x}^{n}}$ of the cotangent bundle as
${\displaystyle g=g_{ij}\,\mathrm {d} x^{i}\otimes \mathrm {d} x^{j}}$.
Endowed with this metric, the differentiable manifold $(M, g)$ is a Riemannian manifold. 

The starting point of Riemannian geometry is Riemannian metric, geodesics can be obtained by variation through Riemannian metric. More specifically,
$$S=\int_{x^{1}}^{x^{2}} dS=\int_{x^{1}}^{x^{2}} \sqrt{g_{ij} d x^{i} d x^{j}}$$So to make $\delta S = 0$, there must be
$$\frac{1}{2} \frac{\partial g_{ij}}{\partial x^{l}} \frac{d x^{i}}{dt} d x^{j}-\frac{\partial g_{il}}{\partial x^{j}} \frac{d x^{i}}{dt} d x^{j}-g_{i l} d\left(\frac{d x^{i}}{dt}\right)=0$$By arranging, the full geodesic equation is $$\frac{d^{2} x^{i}}{d t^{2}}+\Gamma_{lk}^{i} \frac{d x^{l}}{dt} \frac{d x^{k}}{dt}=0$$
We can change the equation of geodesic to
\begin{equation}\label{eq17}
  d\left(\frac{d x^{i}}{dt}\right)=-\Gamma_{lk}^{i} \frac{d x^{l}}{dt} d x^{k}
\end{equation}

\begin{definition}\cite{1,16}
  If for a given smooth non-degenerate symmetric (0,2)-tensor field $g$ everywhere on the $m$ dimensional smooth manifold $M$, then $M$ is called the generalized Riemannian manifold. $g$ is metric tensor on $M$.
\end{definition}
If $g$ is positive definite, then $M$ is Riemannian manifold.

\begin{theorem}\cite{1,15,16}
There must be a Riemann metric on the $m$-dimensional smooth manifold $M$.
\end{theorem}

Let $M$ be a smooth manifold and let $C^{\infty}(M, TM)$ be the space of vector fields on $M$, that is, the space of smooth sections of the tangent bundle $TM$. Then an affine connection on $M$ is a bilinear map
\begin{align}
C^{\infty }(M,TM)\times C^{\infty }(M,TM)&\rightarrow C^{\infty }(M,TM)\notag\\
(X,Y)&\mapsto \nabla _{X}Y\notag
\end{align}
such that for all smooth functions $f$ in $C^{\infty}(M, \mathbb{R})$ and all vector fields $X, Y$ on $M$:
\begin{enumerate}
  \item ${{\nabla}_{fX}}Y=f{{\nabla}_{X}}Y$, that is, $\nabla$ is $C^{\infty}(M, \mathbb{R})$ linear in the first variable;
  \item ${{\nabla}_{X}}\left( fY \right)=df\left( X \right)Y+f{{\nabla}_{X}}Y$, that is, $\nabla$ satisfies Leibniz rule in the second variable.
\end{enumerate}

Accordingly, the curvature tensor is shown as
\begin{align}
  & R:{{C}^{\infty }}\left( TM \right)\times {{C}^{\infty }}\left( TM \right)\times {{C}^{\infty }}\left( TM \right)\to {{C}^{\infty }}\left( TM \right) \notag\\
 & R\left( X,Y \right)Z=\left( {{\nabla }_{X}}{{\nabla }_{Y}}-{{\nabla }_{Y}}{{\nabla }_{X}}-{{\nabla }_{\left[ X,Y \right]}} \right)Z=-R\left( Y,X \right)Z\notag
\end{align}Correspondingly, torsion tensor is given by
\begin{align}
  & T:{{C}^{\infty }}\left( TM \right)\times {{C}^{\infty }}\left( TM \right)\to {{C}^{\infty }}\left( TM \right)\notag \\
 & T\left( X,Y \right)={{\nabla }_{X}}Y-{{\nabla }_{Y}}X-\left[ X,Y \right]=-T\left( Y,X \right)\notag
\end{align}
The local expression of them are given by
\begin{align}
 & T_{ij}^{k}=\Gamma _{ij}^{k}-\Gamma _{ji}^{k}-c_{ij}^{k} \notag\\
 & R_{lij}^{k}=\Gamma _{jl}^{s}\Gamma _{is}^{k}-\Gamma _{il}^{s}\Gamma _{js}^{k}+{{X}_{i}}\Gamma _{jl}^{k}-{{X}_{j}}\Gamma _{il}^{k}-c_{ij}^{s}\Gamma _{sl}^{k}\notag
\end{align}
For natural frame basis,  then $c_{ij}^{k}=0$,  therefore, the formulas above are rewritten as
\begin{align}
  & T_{ij}^{k}=\Gamma _{ij}^{k}-\Gamma _{ji}^{k}\notag \\
 & R_{lij}^{k}=\Gamma _{jl}^{s}\Gamma _{is}^{k}-\Gamma _{il}^{s}\Gamma _{js}^{k}+{{X}_{i}}\Gamma _{jl}^{k}-{{X}_{j}}\Gamma _{il}^{k}\notag
\end{align}

\begin{definition}\cite{1}
  Let $(M, g)$ be a $m$-dimensional generalized Riemannian manifold, $\nabla$ is an affine connection on $M$, if $$\nabla g=0$$then $\nabla$ is admissible connection of generalized Riemannian manifold $(M, g)$.
\end{definition}
Note that the a torsion free connection on a tangent bundle, which holds the Riemann metric constant. The fundamental theorem of Riemann geometry shows that existence of unique ties satisfies these properties. In the theory of Riemann manifold and pseudo Riemann manifold, the expression of the coordinate space of the connection are Christoffel symbol.

\begin{theorem}\cite{1,13}
In a local coordinate system $\left( U,{{x}^{i}} \right)$, connection $\nabla\colon\mathfrak{X}\left( U \right)\times \mathfrak{X}\left( U \right)\to \mathfrak{X}\left( U \right)$ can be arbitrarily determined, law of motion at the natural shelf field $\left\{ \frac{\partial }{\partial {{x}^{i}}};1\le i\le m \right\}$:
\begin{equation}\label{eq16}
  {{\nabla}_{\frac{\partial }{\partial {{x}^{k}}}}}\frac{\partial }{\partial {{x}^{j}}}=\Gamma _{jk}^{i}\frac{\partial }{\partial {{x}^{i}}}
\end{equation}
\end{theorem}

In other words, the connection matrix of connection $\nabla$ is $\omega =\left( \omega _{i}^{k} \right)$ on the local coordinate $x^{i}$,
where
\begin{equation}\label{eq15}
  \omega _{i}^{j}=\Gamma _{ik}^{j}d{{x}^{k}}
\end{equation}
then
\[\nabla g=\left( d{{g}_{ij}}-\omega _{i}^{k}{{g}_{kj}}-\omega _{j}^{k}{{g}_{ki}} \right)\otimes d{{x}^{i}}\otimes d{{x}^{j}}\]
it's equivalent to the equality
\begin{equation}\label{eq7}
  d{{g}_{ij}}=\omega _{i}^{k}{{g}_{kj}}+\omega _{j}^{k}{{g}_{ki}}
\end{equation}
\begin{theorem}\cite{1}
  Let $(M, g)$ be a $m$-dimensional generalized Riemannian manifold, then
  there is a unique affine connection $\nabla$  on $M$ with the following two properties:
  \begin{enumerate}
    \item the connection is torsion-free.
    \item parallel transport is an isometry.
  \end{enumerate}
  This connection is called  Christoffel-Levi-Civita connection.
\begin{proof}
Suppose that $\nabla$ is torsion-free isometric connection on $M$, with \eqref{eq15}, then \[d{{g}_{ij}}=\omega _{i}^{k}{{g}_{kj}}+\omega _{j}^{k}{{g}_{ki}}\]
and $\Gamma _{ik}^{j}=\Gamma _{ki}^{j}$, it leads to the result
\begin{equation}\label{eq4}
  \frac{\partial {{g}_{ij}}}{\partial {{x}^{k}}}={{g}_{lj}}\Gamma _{ki}^{l}+{{g}_{il}}\Gamma _{kj}^{l}
\end{equation}
Rotation index gives
$${\displaystyle {\Gamma ^{l}}_{jk}={\tfrac {1}{2}}g^{lr}\left\{\partial _{k}g_{rj}+\partial _{j}g_{rk}-\partial _{r}g_{jk}\right\}}$$
where as usual $g^{ij}$ are the coefficients of the dual metric tensor, i.e. the entries of the inverse of the matrix ($g_{kl}$).

\end{proof}
\end{theorem}
\begin{definition}\cite{1}
 A smooth curve $\gamma =\left[ a,b \right]\to M$, which satisfies (with  ${\dot{x}^{i}}=\frac{d}{dt}{{x}^{i}}\left( \gamma \left( t \right) \right)$ etc.)
 \[\ddot{x}^{i}(t)+\Gamma_{j k}^{i}(x(t)) \dot{x}^{j}(t) \dot{x}^{k}(t)=0, \text { for } i=1, \ldots, n\]
 is called a geodesic.
\end{definition}
Actually, geodesics are sometimes illustrated as the equilibrium position of a spring on a slippery surface.
\begin{equation}\label{eq13}
  \ddot{x}^{m}+\Gamma_{i j}^{m} \dot{x}^{i} \dot{x}^{j}=0
\end{equation}
where $\dot{x}^{j}=\frac{d{{x}^{j}}}{dt}=v^{j}$ is the velocity.

We call \eqref{eq13} the equation of geodesic. In this equation the properties of the surface appear only through the metric tensor $g$ and its derivatives (via the Christoffel symbol $\Gamma_{i j}^{m}$). This allows us to work in the $m-$dimensional space
with the metric $g$ without any reference to the $n-$dimensional Euclidean space. If the metric tensor $g$ as a function of the point in the space is constant, its derivatives vanish and so do all the Christoffel symbols. The equation of geodesic is then
$$\ddot{x}^{m}=0$$
Obviously, in a flat spacetime similar to Euclidean space, the metric is equal everywhere, so the derivative of the metric to the coordinate is zero. It is easy to see that the geodesic line at this time is a straight line.
In this special case along with the linear solution $x=at+b$, where $a,b$ are constants.

\subsection{Riemannian curvature tensor}
This section is mainly based on \cite{1,14,15}.
Let $M$ be an $n$-dimensional complete Riemannian manifold with the Riemannian metric $g_{ij}$. The Levi-Civita connection is given by the Christoffel symbols
\begin{equation}\label{eq1}
  \Gamma _{ij}^{k}=\frac{1}{2}{{g}^{kl}}\left( \frac{\partial {{g}_{jl}}}{\partial {{x}^{i}}}+\frac{\partial {{g}_{il}}}{\partial {{x}^{j}}}-\frac{\partial {{g}_{ij}}}{\partial {{x}^{l}}} \right)
\end{equation}
where $g^{ ij}$ is the inverse of $g_{ ij}$ . The summation convention of summing over repeated indices is used here and throughout the paper.
The Riemannian curvature tensor is given by
\[R_{ijl}^{k}=\frac{\partial \Gamma _{jl}^{k}}{\partial {{x}^{i}}}-\frac{\partial \Gamma _{il}^{k}}{\partial {{x}^{j}}}+\Gamma _{ip}^{k}\Gamma _{jl}^{p}-\Gamma _{jp}^{k}\Gamma _{il}^{p}\]
We lower the index to the third position, so that $${{R}_{ijkl}}={{g}_{kp}}R_{ijl}^{p}$$
The curvature tensor $R_{ ijkl}$ is anti-symmetric in the pairs $i, j$ and $k, l$ and symmetric in their interchange: $${{R}_{ijkl}}=-{{R}_{jikl}}=-{{R}_{ijlk}}={{R}_{klij}}$$
The Ricci tensor is the contraction
\[
  {{R}_{ik}}=R_{ijk}^{j}={{g}^{jl}}{{R}_{ijkl}}
\]
and mixed (1,1) curvature tensor is
$R_{i}^{j}={{g}^{jk}}{{R}_{ik}}$.
The determinant of $R^i_j$ is invariant under any coordinate transformation. Similarly, the trace $\mathrm{tr} (R^i_j) =R^i_i$ of matrix $R^i_j$ is an invariant. Furthermore,  the scalar curvature $R={{g}^{jk}}{{R}_{jk}}$ which is a ${{C}^{\infty }}\left( M \right)$ scalar function on $M$.
We denote the covariant derivative of a vector field $v=v^{j} \frac{\partial}{\partial x^{j}}$ by
\begin{equation}\label{eq2}
  \nabla_{i} v^{j}=\frac{\partial v^{j}}{\partial x^{i}}+\Gamma_{i k}^{j} v^{k}
\end{equation}
and of a 1-form by
\[\nabla_{i} v_{j}=\frac{\partial v_{j}}{\partial x^{i}}-\Gamma_{i j}^{k} v_{k}\]where $\frac{d{{x}^{k}}}{dt}={{v}^{k}},~{{v}_{k}}={{v}^{l}}{{g}_{kl}}$.
These definitions extend uniquely to tensors so as to preserve the product rule and contractions. For the exchange of two covariant derivatives, we have
\begin{align}\label{eq10}
\nabla_{i} \nabla_{j} v^{l}-\nabla_{j} \nabla_{i} v^{l} &=R_{i j k}^{l} v^{k} \\
\nabla_{i} \nabla_{j} v_{k}-\nabla_{j} \nabla_{i} v_{k} &=R_{i j k l} g^{l m} v_{m}\notag
\end{align}
and similar formulas for more complicated tensors. The second Bianchi identity is
\begin{equation}\label{eq43}
  R_{ijk;l}^{h}+R_{ikl;j}^{h}+R_{ilj;k}^{h}=0
\end{equation}
or \[\nabla_{m} R_{i j k l}+\nabla_{i} R_{j m k l}+\nabla_{j} R_{m i k l}=0\]
where $R_{ijk;l}^{h}={{\nabla}_{l}}R_{ijk}^{h}$.

\subsection{Cartan Structural Equations}

Let $\nabla$ is a connection on the tangent bundle $TM$, in local coordinates, by \eqref{eq16},  we can obtain
$\nabla\frac{\partial }{\partial {{x}^{j}}}=\Gamma _{lj}^{k}d{{x}^{l}}\otimes \frac{\partial }{\partial {{x}^{k}}}$, note that
$\omega _{j}^{k}=\Gamma _{ij}^{k}d{{x}^{i}}$, then \cite{1,14,15,16} \[\nabla\frac{\partial }{\partial {{x}^{j}}}=\omega _{j}^{k}\otimes \frac{\partial }{\partial {{x}^{k}}}\]
for the representation of the connection in a general frame, suppose that ${{e}_{1}},\cdots ,{{e}_{m}}$ is a tangent frame field of $M$ on open set $U$, ${{\omega }^{1}},\cdots ,{{\omega }^{m}}$ is dual frame field of ${{e}_{1}},\cdots ,{{e}_{m}}$, for each point on $U$, there has ${{\omega }^{i}}\left( {{e}_{j}} \right)=\delta _{j}^{i}$. Note that
\[\nabla{{e}_{j}}=\omega _{j}^{k}{{e}_{k}}\]
\[{{\nabla}_{{{e}_{i}}}}{{e}_{j}}=\omega _{j}^{k}\left( {{e}_{i}} \right){{e}_{k}}=\Gamma _{ji}^{k}{{e}_{k}}\]
and $\omega _{j}^{k}=\Gamma _{ij}^{k}{{\omega }^{i}}$ is the connection form of $\nabla$ under the frame field ${{e}_{1}},\cdots ,{{e}_{m}}$.

\begin{theorem}[Cartan Structural Equations]\label{t2}\cite{1,15,16}
Let $n$ dimensional affine connection space $\left( M,\nabla \right)$ be under the local tangent frame field $\left( U,{{e}_{i}} \right)$,dual linear differential form of cotangent frame field, connection form, torsion form, curvature form are respectively shown as
${{\omega }^{i}},\omega _{j}^{i},{{\Omega }^{i}},\Omega _{j}^{i},1\le i,j\le n$, then
geometric structure equation are given by
\begin{align}
  & {{\Omega }^{i}}=d{{\omega }^{i}}-{{\omega }^{j}}\wedge \omega _{j}^{i}=\frac{1}{2}T_{kl}^{i}{{\omega }^{k}}\wedge {{\omega }^{l}} \notag\\
 & d{{\Omega }^{i}}={{\omega }^{j}}\wedge \Omega _{j}^{i}-{{\Omega }^{j}}\wedge \omega _{j}^{i}  \notag\\
 & \Omega _{i}^{j}=d\omega _{i}^{j}-\omega _{i}^{k}\wedge \omega _{k}^{j}=\frac{1}{2}R_{ikl}^{j}{{\omega }^{k}}\wedge {{\omega }^{l}}\notag\\
 & d\Omega _{i}^{j}=\omega _{i}^{k}\wedge \Omega _{k}^{j}-\Omega _{i}^{k}\wedge \omega _{k}^{j}  \notag
\end{align}
\end{theorem}
The Cartan structural equations can be rewritten in a matrix form
\begin{align}
  & \Theta =d\vartheta -\vartheta \wedge \omega  \notag\\
 & d\Theta =\vartheta \wedge \Omega -\Theta \wedge \omega  \notag\\
 & \Omega =d\omega -\omega \wedge \omega  \notag\\
 & d\Omega =\omega \wedge \Omega -\Omega \wedge \omega  \notag
\end{align}
if we denote
$\omega =\left( \omega _{i}^{k} \right),~\vartheta =\left( {{\omega }^{i}} \right),~\Omega =\left( \Omega _{k}^{j} \right),~\Theta =\left( {{\Omega }^{i}} \right)$.

\section{Geospin Matrix and Geospin Variables}
In this section, we present the mathematical formalism involved in \cite{13},    it shows the existence of matrix called geospin matrix for simplifying the equation of geodesics, in the case of Riemannian and pseudo-Riemannian manifolds.

As we know, geodesics are commonly seen in the study of Riemannian geometry and more generally metric geometry. In general relativity, geodesics in spacetime describe the motion of point particles under the influence of gravity alone. In particular, the path taken by a falling rock, an orbiting satellite, or the shape of a planetary orbit are all geodesics in curved spacetime. Obviously, it has wide applications in reality. More generally, the topic of sub-Riemannian geometry deals with the paths that objects may take when they are not free, and their movement is constrained in various ways.

In \cite{13}, a new variable matrix function--geospin matrix $W$ that is defined by Levi-Civita connection along with velocity field. The geodesic equation problem is transformed into a linear dynamic system problem. The geodesic equation is greatly simplified, the dynamic solution can be obtained in some way, and the related problems of the matrix composed of new variables can be dynamically discussed.

In Riemannian geometry, the fundamental theorem of Riemannian geometry states that on any Riemannian manifold or pseudo-Riemannian manifold there is a unique torsion-free metric connection, called the Levi-Civita connection of the given metric. Here a metric or Riemannian connection is a connection which preserves the metric tensor. By using Levi-Civita connection of the given metric, we can define a new variable below. 
\begin{definition}\cite{13}
Let $(M,g)$ be a Riemannian manifold with Eq\eqref{eq1}, two kinds of geospin variables can be defined as
  \begin{equation}
 W_{i}^{j}=\Gamma _{ik}^{j}{{v}^{k}},~~~{{W}_{ik}}=\Gamma _{ik}^{j}{{v}_{j}}
\end{equation}
where  ${{v}^{j}}, {{v}_{j}}$ are the components of  $v={{v}^{j}}\frac{\partial }{\partial {{x}^{j}}}={{v}_{j}}\frac{\partial }{\partial {{x}_{j}}}=\frac{dx}{dt}$ respectively.
\end{definition}
As geospin variables defined above, we mainly study the (1,1) form geospin variable $W_{i}^{j}$ that relates to the geospin matrix $W$ below.

\begin{definition}\cite{13}
The geospin matrix can be defined as
  $$W=\left( W_{j}^{k} \right)$$
where $W_{j}^{k}$ are geospin variables.
\end{definition}
Obviously, the symmetric holds ${{W}_{kj}}={{W}_{jk}}$. The covariant derivative of a vector field $v$ can be rewritten as follows \cite{13}
\begin{equation}\label{eq11}
  {{\nabla }_{k}}{{v}^{j}}=\frac{\partial {{v}^{j}}}{\partial {{x}^{k}}}+W_{k}^{j},~~{{\nabla }_{k}}{{v}_{j}}=\frac{\partial {{v}_{j}}}{\partial {{x}^{k}}}-{{W}_{kj}}
\end{equation}
Namely, Eq \eqref{eq2} has rewritten as simple form. As a consequence, the geodesic for Riemannian manifolds can be rewritten in a compact and simple form by using the geospin variables as follows:
\begin{equation}\label{eq3}
  \left\{ \begin{matrix}
   \frac{d{{x}^{k}}}{dt}={{v}^{k}}\begin{matrix}
   {} &  \\
\end{matrix}  \\
   \frac{d{{v}^{k}}}{dt}+W_{j}^{k}{{v}^{j}}=0 \\
\end{matrix} \right.
\end{equation}
This generalizes the notion of geodesic for Riemannian manifolds. However, in metric geometry the geodesic considered is often equipped with natural parameterization. The local existence and uniqueness theorem for geodesics states that geodesics on a smooth manifold with an affine connection exist, and are unique. More precisely, geodesic equation Eq\eqref{eq3} based on new matrix variable is written in a matrix form as
\begin{equation}
  \left\{ \begin{matrix}
   \frac{d{{x}}}{dt}={{v}}\begin{matrix}
   {} &  \\
\end{matrix}  \notag \\
   \frac{d{{v}}}{dt}=-Wv   \notag\\
\end{matrix} \right.
\end{equation}
This is a dynamic system. This formulation of the geodesic equation of motion can be useful for computer calculations and to compare general relativity with Newtonian gravity.  The left side of the equation is the acceleration vector of the curve on the manifold, so the equation means that the geodesic line is a curve with zero acceleration on the manifold, so the geodesic line must be a constant velocity curve.

The local existence and uniqueness theorem for geodesics states that geodesics on a smooth manifold with an affine connection exist, and are unique.
Hence, the geodesics equation is transformed to initial problem \cite{13}
\begin{equation}\label{eq14}
  \frac{dv}{dt}=-Wv,~~~v\left( {{t}_{0}} \right)={{v}_{0}}
\end{equation}
and $\frac{\partial W}{\partial t}=0$. The proof of this theorem follows from the theory of ordinary differential equations, by noticing that the geodesic equation is a second-order ODE. Existence and uniqueness then follow from the Picard-Lindel\"{o}f theorem for the solutions of ODEs with prescribed initial conditions. Suppose that $W$ is a constant matrix, that is, $W=W_{0}$, then
 it has unique solution formally given by  $v\left( t \right)={{v}_{0}}{{e}^{-W_{0}t}}$. With a new method based on the theory of geospin variables, we solve geodesic equations in a accurate form $v\left( t \right)={{v}_{0}}{{e}^{-W_{0}t}}$,  by Taylor expansion, it has \[v={{v}_{0}}{{e}^{-W_{0}t}}={{v}_{0}}\left( I-W_{0}t+\frac{1}{2}{{W_{0}}^{2}}{{t}^{2}}+\cdots  \right)\]We denote
\[\widehat{u}\left( W_{0},t \right)={{\widehat{u}}_{t}}\left( W_{0} \right)=-{{v}_{0}}\int{W_{0}tdt}+\frac{1}{2}{{v}_{0}}\int{{{W_{0}}^{2}}{{t}^{2}}dt}+\cdots \]as a nonlinear term of the coordinate, and satisfying $\widehat{u}\left( 0,t \right)={{\widehat{u}}_{t}}\left( 0 \right)=0$, where 0 is zero matrix here,  Then the nonlinear coordinate solution of geodesic equation can be written as a form \[u\left( t \right)={{v}_{0}}t+{{u}_{0}}+{{\widehat{u}}_{t}}\left( W_{0} \right)\]Therefore, the nonlinear coordinate solution only depends on the concrete form of geospin matrix $W$. Indeed, in generally, the geospin matrix $W$ is not a constant matrix, obviously, it brings complexity to solve and realize how precisely the solution of geodesics equation \eqref{eq14} would be.

One the other way, let ${{W}_{ij}}^{*}={{g}_{ki}}W_{j}^{k}$ be denoted, then according to \eqref{eq4}, it yields
\[\frac{\partial {{g}_{ij}}}{\partial {{x}^{k}}}{{v}^{k}}={{g}_{lj}}\Gamma _{ki}^{l}{{v}^{k}}+{{g}_{il}}\Gamma _{kj}^{l}{{v}^{k}}={{g}_{lj}}W_{i}^{l}+{{g}_{il}}W_{j}^{l}\]
That is \cite{13},
\begin{equation}\label{eq5}
  \frac{\partial {{g}_{ij}}}{\partial {{x}^{k}}}{{v}^{k}}={{W}_{ji}}^{*}+{{W}_{ij}}^{*}
\end{equation}
the diagonal element of geospin matrix $W$ are expressed in the form
$${{w}^{\left( k \right)}}=W_{k}^{k}=\Gamma _{ki}^{k}{{v}^{i}}=\frac{1}{2}{{g}^{rp}}\frac{\partial {{g}_{rp}}}{\partial {{x}^{i}}}{{v}^{i}}$$
Then by using \eqref{eq5}, the diagonal element can be shown as
\[{{w}^{\left( k \right)}}=W_{k}^{k}=\frac{1}{2}\left( {{g}^{rp}}{{W}_{rp}}^{*}+{{g}^{rp}}{{W}_{pr}}^{*} \right)\]
\eqref{eq17} can be rewritten in a form
\[d{{v}^{i}}=-W_{k}^{l}{{\omega }^{k}}=-{{q}^{i}}dt\]That is,
\[\frac{d{{v}^{i}}}{dt}=-{{q}^{i}}\]

So now, let's study the geospin variables more precisely,
\[W_{i}^{k}=\Gamma _{ij}^{k}{{v}^{j}}=\frac{1}{2}{{g}^{kl}}\left( {{v}^{j}}\frac{\partial {{g}_{jl}}}{\partial {{x}^{i}}}+{{v}^{j}}\frac{\partial {{g}_{il}}}{\partial {{x}^{j}}}-{{v}^{j}}\frac{\partial {{g}_{ij}}}{\partial {{x}^{l}}} \right)\]
Plugging \eqref{eq5} into it, it brings equation to us,
\begin{align}
  W_{i}^{k}& =\frac{1}{2}{{g}^{kl}}\left( W_{il}^{*}+W_{li}^{*}+{{v}^{j}}\frac{\partial {{g}_{jl}}}{\partial {{x}^{i}}}-{{v}^{j}}\frac{\partial {{g}_{ij}}}{\partial {{x}^{l}}} \right) \notag\\
 & =\frac{1}{2}{{g}^{kl}}\left( W_{il}^{*}+W_{li}^{*} \right)+\frac{1}{2}{{g}^{kl}}\left( {{v}^{j}}\frac{\partial {{g}_{jl}}}{\partial {{x}^{i}}}-{{v}^{j}}\frac{\partial {{g}_{ij}}}{\partial {{x}^{l}}} \right) \notag
\end{align}Furthermore, we get
\[W_{i}^{k}{{v}_{k}}=\frac{1}{2}{{v}^{l}}{{v}^{j}}\frac{\partial {{g}_{jl}}}{\partial {{x}^{i}}}\]
where we have used identity
\begin{align}
{{v}^{l}}{{v}^{j}}\frac{\partial {{g}_{jl}}}{\partial {{x}^{i}}}  &={{v}^{l}}{{v}^{j}}{{g}_{pl}}\Gamma _{ij}^{p}+{{v}^{l}}{{v}^{j}}{{g}_{pj}}\Gamma _{il}^{p} \notag\\
 & ={{v}^{l}}{{g}_{pl}}W_{i}^{p}+{{v}^{j}}{{g}_{pj}}W_{i}^{p} \notag\\
 & =2{{v}_{p}}W_{i}^{p} \notag
\end{align}Thus, we obtain
\[2W_{i}^{k}{{v}_{k}}={{v}^{l}}{{v}^{j}}\frac{\partial {{g}_{jl}}}{\partial {{x}^{i}}}=2{{v}_{p}}W_{i}^{p}\]Going further leads to the identity
\[W_{i}^{k}{{v}_{k}}{{v}^{i}}=\frac{1}{2}{{v}^{l}}{{v}^{j}}{{v}^{i}}\frac{\partial {{g}_{jl}}}{\partial {{x}^{i}}}=\frac{1}{2}{{v}^{l}}{{v}^{j}}\left( W_{jl}^{*}+W_{lj}^{*} \right)\]or the form
\[2W_{i}^{k}{{v}_{k}}{{v}^{i}}={{v}^{l}}{{v}^{j}}\left( W_{jl}^{*}+W_{lj}^{*} \right)\]

\section{Geometrodynamics on Riemannian manifolds}
In this section, we will discuss the geometrodynamics on Riemannian manifolds $M$. By analyzing the geodesic, it surely implies that geospin matrix $W=\left( W_{j}^{k} \right)$ is a dynamic system alone. Obviously, it's a truly geometrodynamics on Riemannian manifolds.

Based on Eq\eqref{eq3}, we denote ${{q}^{k}}=W_{i}^{k}{{v}^{i}}$ as geometric acceleration , then the geodesic equation can be simply written in the form
\[{{q}^{k}}=-\frac{d{{v}^{k}}}{dt}\]
Furthermore, we can construct an invariance:
$Q={{q}^{k}}{{v}_{k}}=W_{i}^{k}{{v}^{i}}{{v}_{k}}$.

Geometrodynamics on Riemannian manifolds should be built on the geometric structure equation as theorem \ref{t2} given.  Essentially, we can see that all geometric structural equations are related to the basic quantities ${{\omega }^{k}},\omega _{i}^{j}$ which can be expressed in a specific form
\begin{align}\label{eq6}
  & {{\omega }^{i}}=d{{x}^{i}} \\
 & \omega _{i}^{j}=\Gamma _{ik}^{j}d{{x}^{k}}=\Gamma _{ik}^{j}{{\omega }^{k}}\notag
\end{align}
By using the dynamical expression, that is velocity field and the geospin variable given by  $\left\{ {{v}^{i}}=\frac{d{{x}^{i}}}{dt},~~W_{i}^{j}=\Gamma _{ik}^{j}{{v}^{k}} \right\}$, \eqref{eq6} can be rewritten as
\begin{align}\label{eq8}
  & {{\omega }^{i}}={{v}^{i}}dt \\
 & \omega _{i}^{j}=W_{i}^{j}dt \notag
\end{align}
\begin{definition}\label{d1}
 The quantities ${{v}^{i}},~W_{i}^{j}$ are defined as dynamical variables.
\end{definition}
Helping with \eqref{eq8} leads to reconsider the geometric structural equations \ref{t2} that are given
\begin{align}\label{eq9}
  & {{\Omega }^{i}}=d{{\omega }^{i}}-{{\omega }^{j}}\wedge \omega _{j}^{i}=\frac{1}{2}T_{kl}^{i}{{\omega }^{k}}\wedge {{\omega }^{l}}=\frac{1}{2}T_{kl}^{i}{{v}^{k}}\wedge {{v}^{l}}{{\left( dt \right)}^{2}}, \\
 & d{{\Omega }^{i}}={{\omega }^{j}}\wedge \Omega _{j}^{i}-{{\Omega }^{j}}\wedge \omega _{j}^{i}  \notag\\
 & \Omega _{i}^{j}=d\omega _{i}^{j}-\omega _{i}^{k}\wedge \omega _{k}^{j}=\frac{1}{2}R_{ikl}^{j}{{\omega }^{k}}\wedge {{\omega }^{l}}=\frac{1}{2}R_{ikl}^{j}{{v}^{k}}\wedge {{v}^{l}}{{\left( dt \right)}^{2}}  \notag\\
 & d\Omega _{i}^{j}=\omega _{i}^{k}\wedge \Omega _{k}^{j}-\Omega _{i}^{k}\wedge \omega _{k}^{j}  \notag
\end{align}
Using \eqref{eq8} to give a dynamical interpretation of the geometrodynamics on Riemannian manifolds is based on the geometric structural equations \eqref{eq9} above. To write \eqref{eq8} in a matrix form is given by
\[\left( \begin{matrix}
   {{\omega }^{i}}  \\
   \omega _{i}^{j}  \\
\end{matrix} \right)=\left( \begin{matrix}
   {{v}^{i}}  \\
   W_{i}^{j}  \\
\end{matrix} \right)dt\]
As a consequence, we will mainly focus on definition \ref{d1} to rewrite pure geometric equation which leads to the geometrodynamics structural equations.

\subsection{Geometrodynamics structural equations}
\begin{theorem}[Geometrodynamics structural equations]\label{t3}
The geometrodynamics structural equations for $\left\{ {{\Omega }^{i}},\Omega _{i}^{j} \right\}$ are \[\alpha _{i}^{j}-W_{i}^{k}\wedge W_{k}^{j}=\frac{1}{2}R_{jkl}^{i}{{v}^{k}}\wedge {{v}^{l}}\]
and
\[{{a}^{i}}-{{v}^{j}}\wedge W_{j}^{i}=\frac{1}{2}T_{jk}^{i}{{v}^{j}}\wedge {{v}^{k}}\]
where $\alpha _{i}^{j}=d{W_{i}^{j}}/dt$ and ${{a}^{i}}=d{{v}^{i}}/dt$.

The geometrodynamics structural equations for $\left\{ d{{\Omega }^{i}},d\Omega _{i}^{j} \right\}$ are given by
\begin{align}\label{eq18}
  & \frac{d{{\Omega }^{i}}}{dt}={{v}^{j}}\wedge \Omega _{j}^{i}-{{\Omega }^{j}}\wedge W_{j}^{i} \\
 & \frac{d\Omega _{i}^{j}}{dt}=W_{i}^{k}\wedge \Omega _{k}^{j}-\Omega _{i}^{k}\wedge W_{k}^{j} \notag
\end{align}
\begin{proof}
 Plugging \eqref{eq8} into the geometric structural equations \eqref{eq9}, we can easily obtain the geometrodynamics structural equations on Riemannian manifolds.

\end{proof}
\end{theorem}

\begin{corollary}\label{c1}
If ${{\Omega }^{i}}=0,~\Omega _{i}^{j}=0$ hold, then
  \[\alpha _{i}^{j}=W_{i}^{k}\wedge W_{k}^{j}\]
  \[{{a}^{i}}={{v}^{j}}\wedge W_{j}^{i}\]
\begin{proof}
  If ${{\Omega }^{i}}=0,~\Omega _{i}^{j}=0$ hold,  by directly calculation, we have
\[d{{\omega }^{i}}={{\omega }^{j}}\wedge \omega _{j}^{i}={{a}^{i}}dtdt={{v}^{j}}\wedge W_{j}^{i}dtdt\]
\[d\omega _{i}^{j}=\omega _{i}^{k}\wedge \omega _{k}^{j}=\frac{d W_{i}^{j}}{d t}dtdt=W_{i}^{k}\wedge W_{k}^{j}dtdt\]

\end{proof}
\end{corollary}
Note that the second equation in corollary 2 is rightly describing the  geodesic equation, it implies that the geodesic equation is a equation without considering the torsion form that is compatible with  torsion-free.

\begin{corollary}
  If we set $\alpha =\left( \alpha _{k}^{j} \right), a=\left( {{a}^{i}} \right),v=\left( {{v}^{j}} \right)$,  the geometrodynamics structural equations \ref{t3} can also be rewritten in a matrix form
\begin{align}
  & \Omega /dt^2=\alpha -W\wedge W \notag\\
 & \Theta /dt^2 =a-v\wedge W\notag\\
 & d\Theta /dt=v\wedge \Omega -\Theta \wedge W \notag\\
 & d\Omega /dt=W\wedge \Omega -\Omega \wedge W \notag
\end{align}
Subsequently, the corollary \ref{c1} is written as
\begin{align}
  & \alpha =W\wedge W \notag\\
 & a=v\wedge W \notag
\end{align}
\end{corollary}
Actually, the specific details of \eqref{eq18} can be obtained by substituting
geometric structural equations \eqref{eq9} into \eqref{eq18},
\begin{align}
  \frac{d\Omega _{i}^{j}}{dt}& =W_{i}^{k}\wedge \Omega _{k}^{j}-\Omega _{i}^{k}\wedge W_{k}^{j} \notag\\
 & =W_{i}^{k}\wedge \left( d\omega _{k}^{j}-\omega _{k}^{p}\wedge \omega _{p}^{j} \right)-\left( d\omega _{i}^{k}-\omega _{i}^{q}\wedge \omega _{q}^{k} \right)\wedge W_{k}^{j} \notag\\
 & =W_{i}^{k}\wedge d\omega _{k}^{j}-W_{i}^{k}\wedge \omega _{k}^{p}\wedge \omega _{p}^{j}-d\omega _{i}^{k}\wedge W_{k}^{j}+\omega _{i}^{q}\wedge \omega _{q}^{k}\wedge W_{k}^{j} \notag\\
 & =\left( W_{i}^{k}\wedge \alpha _{k}^{j}-\alpha _{i}^{k}\wedge W_{k}^{j} \right)dt^2+\left( W_{i}^{q}\wedge W_{q}^{k}\wedge W_{k}^{j}-W_{i}^{k}\wedge W_{k}^{p}\wedge W_{p}^{j} \right)dt^2 \notag\\
 & =\left( W_{i}^{k}\wedge \alpha _{k}^{j}-\alpha _{i}^{k}\wedge W_{k}^{j} \right)dt^2 \notag
\end{align}
Then it follows \[\frac{d\Omega _{i}^{j}}{d{{t}^{3}}}=W_{i}^{k}\wedge \alpha _{k}^{j}-\alpha _{i}^{k}\wedge W_{k}^{j}\]
In a matrix form, it's given by
\[\frac{d\Omega }{d{{t}^{3}}}=W\wedge \alpha -\alpha \wedge W\]
Similarly, we have
\begin{align}
 \frac{d{{\Omega }^{i}}}{dt} &={{v}^{j}}\wedge \Omega _{j}^{i}-{{\Omega }^{j}}\wedge W_{j}^{i} \notag\\
 & ={{v}^{j}}\wedge d\omega _{i}^{j}-{{v}^{j}}\wedge \omega _{j}^{q}\wedge \omega _{q}^{i}-\left( d{{\omega }^{j}}-{{\omega }^{p}}\wedge \omega _{p}^{j} \right)\wedge W_{j}^{i} \notag\\
 & =\left( {{v}^{j}}\wedge \alpha _{i}^{j}-{{a}^{j}}\wedge W_{j}^{i} \right)dt^2 \notag
\end{align}
It derives \[\frac{d{{\Omega }^{i}}}{d{{t}^{3}}}={{v}^{j}}\wedge \alpha _{i}^{j}-{{a}^{j}}\wedge W_{j}^{i}\]or in a matrix form
\[\frac{d\Theta }{d{{t}^{3}}}=v\wedge \alpha -a\wedge W\]
where we have used the replacement $d\omega _{k}^{j}=\alpha _{k}^{j}dt^2,~~d{{\omega }^{j}}={{a}^{j}}dt^2$, more precisely, \[d\omega _{k}^{j}=d\left( W_{k}^{j}dt \right)=\alpha _{k}^{j}d{{t}^{2}}+W_{k}^{j}{{d}^{2}}t=\alpha _{k}^{j}d{{t}^{2}}\]
In conclusions,
\begin{align}
  & \frac{d\Omega _{i}^{j}}{d{{t}^{3}}}=W_{i}^{k}\wedge \alpha _{k}^{j}-\alpha _{i}^{k}\wedge W_{k}^{j} \notag\\
 & \frac{d{{\Omega }^{i}}}{d{{t}^{3}}}={{v}^{j}}\wedge \alpha _{i}^{j}-{{a}^{j}}\wedge W_{j}^{i} \notag
\end{align}Matrix form is
\begin{align}
  & \frac{d\Omega }{d{{t}^{3}}}=W\wedge \alpha -\alpha \wedge W \notag\\
 & \frac{d\Theta }{d{{t}^{3}}}=v\wedge \alpha -a\wedge W \notag
\end{align}
Above all,
\begin{corollary}\label{c2}
  The geometrodynamics structural equations \ref{t3} can show as
\begin{align}
  & \Theta /d{{t}^{2}}=a-v\wedge W \notag\\
 & d\Theta /d{{t}^{3}}=v\wedge \alpha -a\wedge W\notag\\
 & \Omega /d{{t}^{2}}=\alpha -W\wedge W \notag\\
 & d\Omega /d{{t}^{3}}=W\wedge \alpha -\alpha \wedge W\notag
\end{align}
\end{corollary}
We have to say that $\left\{ a,\alpha ,v,W \right\}$ are dynamical variables.
More precisely, $a,\alpha ,v,W $ represent the acceleration, geometric angular acceleration, velocity, geospin matrix.

The geometrodynamics on Riemannian is definitely built these four dynamical variables, indeed. The corollary \ref{c2} reveals the dynamical essence of Cartan structural equations.  The corollary \ref{c2} can be formally rewritten in a matrix form that is convenient to memorize.
\begin{align}
  & \left( \begin{matrix}
   \Omega   \\
   \Theta   \\
\end{matrix} \right)/d{{t}^{2}}=\left( \begin{matrix}
   \alpha   \\
   a  \\
\end{matrix} \right)-\left( \begin{matrix}
   W  \\
   v  \\
\end{matrix} \right)\wedge W \notag\\
 & \left( \begin{matrix}
   d\Omega   \\
   d\Theta   \\
\end{matrix} \right)/d{{t}^{3}}=\left( \begin{matrix}
   W  \\
   v  \\
\end{matrix} \right)\wedge \alpha -\left( \begin{matrix}
   \alpha   \\
   a  \\
\end{matrix} \right)\wedge W \notag
\end{align}
If let \[\left( \begin{matrix}
   \Omega   \\
   \Theta   \\
\end{matrix} \right)/d{{t}^{2}}=0\] be given, then the corollary \ref{c1} reappears
\[\left( \begin{matrix}
   \alpha   \\
   a  \\
\end{matrix} \right)=\left( \begin{matrix}
   W  \\
   v  \\
\end{matrix} \right)\wedge W\]
To use \eqref{eq10} to consider \eqref{eq11}, then
\begin{align}
  & {{\nabla }_{i}}\left( {{\partial }_{j}}{{v}^{k}}+W_{j}^{k} \right)-{{\nabla }_{j}}\left( {{\partial }_{i}}{{v}^{k}}+W_{i}^{k} \right)=R_{ijl}^{k}{{v}^{l}}  \notag\\
 & {{\nabla }_{i}}\left( {{\partial }_{j}}{{v}_{k}}-{{W}_{jk}} \right)-{{\nabla }_{j}}\left( {{\partial }_{i}}{{v}_{k}}-{{W}_{ik}} \right)={{R}_{ijkl}}{{g}^{lm}}{{v}_{m}}={{R}_{ijkl}}{{v}^{l}} \notag
\end{align}
where \[R_{ijl}^{k}{{v}^{l}}={{v}^{l}}{{\partial }_{i}}\Gamma _{jl}^{k}-{{v}^{l}}{{\partial }_{j}}\Gamma _{il}^{k}+\Gamma _{ip}^{k}W_{j}^{p}-\Gamma _{jp}^{k}W_{i}^{p}\]
The above results show that general relativity is no longer the synonym of Riemannian geometry in physics. Even if there is no general relativity, there is Riemannian geometry in physics. The geometrization of mechanics helps us to connect mechanics, field theory and geometry. Riemannian geometry is actually a research framework of geometry. As long as it can be transformed correspondingly, many conclusions of Riemannian geometry can be directly applied, and more abundant and comprehensive contents may be derived.


\end{document}